\documentclass{article}
\usepackage{amssymb,amsmath,amsfonts,epsfig,latexsym}

\newtheorem{theorem}{Theorem}

\newtheorem{lemma}{Lemma}

\newtheorem{remark}{Remark}

\def\Z{\ensuremath{\mathbb{Z}}}
\def\Q{\ensuremath{\mathbb{Q}}} 
\def\P{\ensuremath{\mathbb{P}}}
\def\C{\ensuremath{\mathbb{C}}}
\def\A{\ensuremath{\mathbb{A}}}
\def\R{\ensuremath{\mathbb{R}}}

\def\cA{\ensuremath{\mathcal{A}}}

\def\M{\mathfrak{M}}

\def\O{\ensuremath{\mathcal{O}}}

\def\div{\mathrm{div}}
\def\Hom{\mathrm{Hom}}
\def\int{\mathrm{int}}

\def\Star{\mathrm{Star\,}}
\def\Sym{\mathrm{Sym\,}}

\def\<{\ensuremath{\langle}}
\def\>{\ensuremath{\rangle}}

\def\PP{P\!P}

\begin{document}
    
\title{Equivariant Chow cohomology of toric varieties}           
    
\author{Sam Payne}

\date{}

\maketitle

\begin{abstract}
   We show that the equivariant Chow cohomology ring of a toric
   variety is naturally isomorphic to the ring of integral piecewise
   polynomial functions on the associated fan.  This gives a large class of singular spaces for which localization holds in equivariant Chow cohomology with integer coefficients.  We also compute the equivariant Chow cohomology of toric prevarieties and general complex hypertoric varieties in terms of piecewise polynomial functions.
\end{abstract}

If $X = X(\Delta)$ is a smooth, complete complex toric variety then
the following rings are canonically isomorphic:\ the equivariant singular
cohomology ring $H^{*}_{T}(X)$, the equivariant Chow cohomology ring
$A^{*}_{T}(X)$, the Stanley-Riesner ring $SR(\Delta)$, and the ring of
integral piecewise polynomial functions $\PP^{*}(\Delta)$.  If $X$ is
simplicial but not smooth then $H^{*}_{T}(X)$ may have torsion and the
natural map from $SR(\Delta)$ takes monomial generators to piecewise
linear functions with rational, but not necessarily integral,
coefficients.  In such cases, these rings are not isomorphic, but they
become isomorphic after tensoring with $\Q$.  When $X$ is not
simplicial, there are still natural maps between these rings, for
instance from $A^{*}_{T}(X)_{\Q}$ to $H^{*}_{T}(X)_{\Q}$ and
from $H^{*}_{T}(X)$ to $\PP^{*}(\Delta)$, but these maps are far
from being isomorphisms in general.

The main purpose of this note is to construct a natural isomorphism from $A^{*}_{T}(X)$ to $\PP^{*}(\Delta)$ for an arbitrary toric variety; the map is obtained by restricting a Chow cohomology class to each of the $T$-orbits $O_{\sigma} \subset X$ for cones $\sigma \in \Delta$.  The equivariant Chow cohomology of $O_{\sigma}$ is naturally isomorphic to the ring $\Sym M_\sigma$ of integral polynomial functions on $\sigma$, where $M_\sigma = M / (\sigma^\perp \cap M)$ (for $u \in M$, the image of $u$ in $M_\sigma$ is identified with the first equivariant Chern class of the equivariant line bundle
$\O_{X}(\div \chi^{u})|_{O_{\sigma}}$ in $A^1_T(O_\sigma))$. The ring of integral piecewise
polynomial functions on $\Delta$ is defined by
\[
  \PP^{*}(\Delta) = \{ f : |\Delta| \rightarrow \R : f|_{\sigma} \in
  \Sym M_\sigma \mbox{ for each } \sigma \in \Delta \}.
\]
The map $f \mapsto (f|_{\sigma})_{\sigma \in \Delta}$ identifies
$\PP^{*}(\Delta)$ with a subring of $\bigoplus_{\sigma \in \Delta}
\Sym M_\sigma$:
\[
\PP^*(\Delta) \cong \{ (f_{\sigma})_{\sigma \in \Delta} : f_{\tau} = f_{\sigma}|_{\tau} \mbox{ for } \tau \prec \sigma \}.
\]
We write $\iota_{\sigma}$ for the inclusion of $O_{\sigma}$ in $X$.

\begin{theorem} \label{main}
    Let $X = X(\Delta)$ be a toric variety.  Then $\bigoplus_{\sigma
    \in \Delta} \iota_{\sigma}^{*}$ maps $A^{*}_{T}(X)$ isomorphically
    onto $\PP^{*}(\Delta)$.
\end{theorem}

\noindent It follows that the presheaf $U \mapsto A^*_T(U)$ for $T$-invariant open sets $U \subset X$ is a sheaf in the $T$-equivariant Zariski topology; in the terminology of \cite{BBFK}, $A^*_T$ is a ``sheaf on the fan".

When $X$ is complete, Theorem \ref{main} can be seen as a Chow cohomology version of Goresky-Kottwitz-MacPherson localization in singular cohomology \cite{GKM}: the equivariant Chow cohomology of $X$ injects into the equivariant cohomology of the fixed points, and the image has an explicit description in terms of the stabilizers of the 1-dimensional orbits.  Let $X$ be a complete toric variety with $T$-fixed points $\{x_1, \ldots, x_r\}$ and 1-dimensional $T$-orbits $E_1, \ldots, E_\ell$.  For $j = 1, \ldots, \ell$, let $T_j$ be the stabilizer of any point in $E_j$, let $M_j$ be the character lattice of $T_j$, and let $N_j \subset N$ be the dual lattice of $M_j$.  Define
\[
\beta_j : \bigoplus_{i=1}^r \Sym M \rightarrow \Sym M_j
\]
to be the map given by
\[
 (f_1, \ldots, f_r) \mapsto f_{a_j} |_{N_{a_j}} - f_{b_j} |_{N_{b_j}},
 \]
 where $x_{a_j}$ and $x_{b_j}$ are the two points in the boundary of $E_j$. (The ordering of $a_j$ and $b_j$ is arbitrary; changing the ordering changes $\beta_j$ by a sign, but does not affect the kernel.)  The following version of Theorem \ref{main} for complete toric varieties is analogous to \cite[Theorem 7.2]{GKM}.
 
\begin{theorem}\label{GKM}
Let $X = X(\Delta)$ be a complete $n$-dimensional toric variety.  The restriction map
\[
A^*_T(X) \rightarrow A^*_T(X^T) = \bigoplus_{i=1}^r \Sym M
\]
is injective, and its image is the intersection of the kernels,
\[
A^*_T(X) \cong \bigcap_{j=1}^\ell \mathrm{ker}(\beta_j).
\]
\end{theorem}

\noindent For localization theorems for the equivariant Chow homology groups $A_*^T(X)$, see \cite{Brion2} and \cite{EG2}.

The relationships between piecewise-polynomial functions and various equivariant cohomology theories on toric varieties have been studied extensively, usually with rational coefficients.   For smooth toric varieties, Bifet, De Concini, and Procesi showed that $H^*_T(X)$ is isomorphic to the Stanley-Reisner ring $SR(\Delta)$ \cite[Theorem 8]{BDP} (the connection between Stanley-Reisner rings and piecewise polynomials had already been made by Billera \cite{Billera}).  Brion rediscovered this result with rational coefficients and made the connection with piecewise-polynomial functions explicit \cite[Proposition 2.2]{Brion1}.  Brion and Vergne then extended the isomorphism $H^*_T(X)_\Q \cong \PP^*(\Delta)_\Q$ to simplicial toric varieties \cite[Proposition 3.2]{BV}, and  Brion noted that a similar argument in this case shows that $A^*_T(X)_\Q$ is isomorphic to $\PP^*(\Delta)_\Q$ \cite[Remark 5.4]{Brion2}.  Barthel, Brasselet, Fieseler, and Kaup gave an example of a nonsimplicial toric variety $X$ with nonvanishing $H^3_T(X)_\Q$, showing that there is, in general, no natural isomorphism from $H^*_T(X)_\Q$ to $\PP^*(\Delta)_\Q$, nor to $A^*_T(X)_\Q$ \cite[Example 1.5]{BBFK}.  In the same paper, they showed that the rational equivariant intersection homology $IH^*_T(X)_\Q$ is a module over $\PP^*(\Delta)_\Q$, and that rational equivariant intersection homology is the unique ``minimal extension sheaf" of $\PP^*_\Q$ modules.   Braden and Lunts make essential use of the sheaf of complex-valued piecewise polynomial functions in their new work on toric equivariant Koszul duality \cite{BL} (although they mistakenly assert in Remark 3.6.1 that the global sections of this sheaf compute equivariant singular cohomology).  Brylinski and Zhang have given an elegant and direct proof of the isomorphism $H^*_T(X)_\Q \cong \PP^*(\Delta)_\Q$ for simplicial toric varieties, using a Mayer-Vietoris spectral sequence \cite[Theorem 6.4]{BZ}.  When $X$ is smooth, their argument also gives $H^*_T(X) \cong \PP^*(\Delta)$ with integer coefficients.  Another simple proof of this isomorphism when $X$ is smooth and complete, using only elementary algebra, is due to Fulton and Musta\c{t}\v{a} \cite{TV}.

The main new contribution of this note is the extension of Brion's isomorphism $A^*_T(X)_\Q \cong \PP^*(\Delta)_\Q$ to the nonsimplicial case, and the improvement from rational coefficients to integer coefficients.  The primary tool that makes this possible is Kimura's inductive method for computing Chow rings of singular varieties using envelopes and resolutions of singularities \cite{Kimura}.  Kimura's method was adapted to the equivariant setting by Edidin and Graham \cite{EG1}.  

In the final two sections, we apply the techniques of the main part of this paper to toric prevarieties and complex hypertoric varieties.  We show that the equivariant Chow cohomology of a toric prevariety is naturally isomorphic to the ring of integral piecewise polynomial functions on the associated ``multifan".  Using the observation of Proudfoot and Webster that general hypertoric varieties are affine bundles over toric prevarieties \cite[Remark 3.6]{PW}, we compute the equivariant Chow cohomology rings of these hypertoric varieties, showing that they are isomorphic to those of the underlying toric prevarieties.

\section{Preliminaries}

The Chow cohomology that we use here is the ``operational" theory developed by Fulton and MacPherson \cite{FM} \cite[Chapter 17]{IT}.   

Roughly speaking, the equivariant Chow cohomology $A^*_G(X)$, for a group $G$ and a $G$-space $X$, should be defined as in the Borel construction of equivariant singular cohomology, by taking the Chow cohomology of $X \times_G EG$, where $EG$ is a contractible space on which $G$ acts freely.  Unfortunately, $EG$ does not exist as a finite dimensional algebraic variety in any reasonable sense.  However, Totaro observed that one can construct an arbitrarily good ``approximation" to $EG$ by taking an open subset $U$ of a representation of $G$ such that $G$ acts freely on $U$ and the codimension of the complement of $U$ is sufficiently large \cite{Classifying Spaces}.  Using Totaro's construction, Edidin and Graham systematically developed equivariant intersection theory, defining $A^i_G(X)$ to be the ordinary Chow cohomology group $A^i(X \times_G U)$, where the codimension of the complement of $U$ is greater than $i$, and showing that, up to canonical isomorphisms, this is independent of the choice of $U$ and respects the multiplicative structure.  That is, the multiplication map in equivariant Chow cohomology
\[
A^i_G(X) \otimes A^j_G(X) \rightarrow A^{i+j}_G(X)
\]
 is given by the multiplication map
\[
A^i(X \times_G U) \otimes A^j(X \times_G U) \rightarrow A^{i + j}(X \times_G U),
\]
when the codimension of the complement of $U$ is greater than $i + j$ \cite{EG1}.  

We follow the convention of Edidin and Graham, writing $X_G$ for $X \times_G U$, where the codimension of the complement of $U$ is assumed to be sufficiently large for the given context.

In this paper, we will only consider equivariant cohomology for the action of a torus $T$.  As observed by Brion \cite[Section 2.2]{Brion1}, the approximations $U$ can be chosen such that $U$ has a faithful action of a large torus $T_U$ with a dense orbit, and such that the action of $T$ on $U$ is given by an inclusion $T \hookrightarrow T_U$.  Then $X_T$ has an action of $\tilde{T} := T_U / T$.  In the case where $X$ is a toric variety, this gives $X_T$ the structure of a toric variety with dense torus $\tilde{T}$, and $X_T$ is smooth (resp.\ complete, simplicial) if and only if $X$ is.  It follows that standard results on the ordinary Chow cohomology of toric varieties can be carried over to the equivariant setting immediately.  For instance, if $X$ is a complete toric variety, then $A^i(X)$ is isomorphic to $\Hom(A_i(X), \Z)$ \cite[Proposition 2.4]{FS} \cite[Theorem 2]{Totaro}.  Therefore, $A^i_T(X)$ is isomorphic to $\Hom(A_i(X_T), \Z)$.  In particular, it follows that the equivariant Chow cohomology of a complete toric variety is torsion free.  Similarly, when $X$ is a complete complex toric variety, Totaro has constructed a natural split injection $A^*(X)_\Q \hookrightarrow H^*(X)_\Q$ \cite[Theorem 6]{Totaro}, so it follows that there is also a natural split injection in the equivariant setting $A^*_T(X)_\Q \hookrightarrow H^*_T(X)_\Q$.

As mentioned in the introduction, $H^*_T(X)$ may have torsion when $X$ is a complete simplicial complex toric variety.  This possibility is probably well-known to experts, but we know of no explicit examples in the literature.  The following is an example of a complete toric surface with 2-torsion in $H^3_T$.

\vspace{10 pt}

\noindent \textbf{Example} \ Let $N = \Z^2$, let $\Delta$ be the complete fan in $N_\R$ whose rays are generated by $(\pm1, \pm 1)$, and  let $X = X(\Delta)$ ($X$ is isomorphic to a quotient of $\P^1 \times \P^1$ by $\Z/2 \oplus \Z/2$).  In the Mayer-Vietoris spectral sequence
\[
E_1^{p,q} = \bigoplus_{\sigma_0, \ldots, \sigma_p} H^q_T (U_{\sigma_0} \cap \cdots \cap U_{\sigma_p}) \Rightarrow H^{p+q}_T(X)
\]
considered in \cite[Section~6]{BZ}, $E_2^{0,3}$ and all of the $E_2$ terms southeast of $E_2^{1,2}$ vanish, so $H^3_T(X)$ is isomorphic to $E_2^{1,2}$.  The complex $0 \rightarrow E_1^{0,2} \rightarrow E_1^{1,2} \rightarrow E_1^{2,2}$ looks like
\[
  0 \rightarrow \Z^2 \oplus \Z^2 \oplus \Z^2 \oplus \Z^2 \stackrel{d}{\longrightarrow} \Z \oplus \Z \oplus \Z \oplus \Z \rightarrow 0,
\]  
and one checks easily that if $(n_1, \ldots, n_4)$ is in the image of $d$, then $n_1 + \cdots + n_4$ is even.  Hence $H^3_T(X)$, which is isomorphic to the cokernel of $d$, has nonvanishing 2-torsion.

\section{ Proof of Theorems \ref{main} and \ref{GKM}}

When $X$ is smooth, the isomorphism in Theorem 1 is well-known to experts.  In this case, by \cite[Proposition 4]{EG1}, capping with the equivariant fundamental class gives an isomorphism 
\[
A^i_T(X)  \stackrel{\cap [X]_T} {\longrightarrow} A^T_{n-i}(X).  
\]
(Note that the grading is such that $A^i_T(X)$ vanishes for $i < 0$ and $A^T_j(X)$ vanishes for $j > n$.)  The isomorphism with $\PP^i(\Delta)$ then follows from \cite[Theorem 5.4]{Brion2}. 

We now deduce the general case of Theorem \ref{main} from the smooth
case, using Kimura's inductive methods for computing Chow groups via
envelopes and resolutions of singularities.  This procedure is
relatively simple for toric varieties, since every toric resolution of
singularities is an envelope.  Recall that an envelope is a proper
morphism $\tilde{X} \rightarrow X$ such that every closed subvariety
of $X$ is the birational image of some closed subvariety of
$\tilde{X}$.

\begin{lemma}  Every proper birational toric morphism is an envelope.
\end{lemma}

\noindent \emph{Proof:} If $\pi: \tilde{X} \rightarrow X$ is a proper 
birational toric morphism, with $\tilde{X} = X(\tilde{\Delta})$ and $X =
X(\Delta)$, then $\tilde{\Delta}$ is a subdivision of $\Delta$.  
For any closed subvariety $Y \subset X$, there is a smallest $T$-orbit
$O_{\sigma}$ such that $Y \cap O_{\sigma}$ is dense in $Y$ (the
closure of $O_{\sigma}$ is the intersection of the $T$-invariant
subvarieties of $X$ containing $Y$).  Then for any
maximal cone $\tilde{\sigma}$ in the induced subdivision of $\sigma$, $\pi$ maps $O_{\tilde{\sigma}}$ isomorphically onto
$O_{\sigma}$.  Hence the closure in $\tilde{X}$ of $O_{\tilde{\sigma}}
\cap \pi^{-1}(Y)$ maps birationally onto $Y$.  \hfill $\Box$

\vspace{10 pt}

For any equivariant envelope $\pi: \tilde{X} \rightarrow X$, the induced maps $\pi_T : \tilde{X}_T \rightarrow X_T$ are envelopes \cite[Lemma 3]{EG1}.  Hence the pullback map $\pi^{*}: A^{*}_{T}(X) \rightarrow A^{*}_{T}(\tilde{X})$ is an injection, and the sequence
\[
0 \rightarrow A^*_T(X) \stackrel{\pi^*}{\longrightarrow} A^*_T(\tilde X) \stackrel{p_1^* - p_2^*}{\longrightarrow} A^*_T(\tilde X \times_X \tilde X)
\]
is exact \cite[Theorem 2.3]{Kimura}.  When $\pi$ is birational and an
isomorphism over an open set $U \subset X$, Kimura gives another description of the image of $\pi^*$ which is also useful.

\begin{lemma} \label{equivariant Kimura}
  Let $S_{i}$ be the irreducible
components of $X \smallsetminus U$, $E_{i} = \pi^{-1}(S_i)$, and $\pi_{i}:
E_{i} \rightarrow S_{i}$ the restriction of $\pi$.  Then a class
$\tilde{c} \in A^{*}_{T}(\tilde{X})$ is in the image of $\pi^{*}$ if
and only if $\tilde{c}|_{E_{i}}$ is in the image of $\pi_{i}^{*}$ for
all $i$.
\end{lemma}

\noindent \emph{Proof:}  Apply
\cite[Theorem 3.1]{Kimura} to the maps $\pi_T: \tilde{X}_T \rightarrow X_T$. \hfill $\Box$

\vspace{10 pt}

Since $E_{i}$ and $S_{i}$ have smaller dimension than $X$,
we can use Lemma \ref{equivariant Kimura} to compute $A^{*}_{T}(X)$ using a resolution of singularities and induction on dimension. 

\vspace{10 pt}

\noindent \emph{Proof of Theorem \ref{main}:}  As observed at the beginning of this section, Theorem \ref{main} is true if $X$ is smooth.  If $X$ is singular, then there is a sequence
\[
  X_r \rightarrow X_{r-1} \rightarrow \cdots \rightarrow X_1  \stackrel{\pi}{\rightarrow} X_0 = X
\]
where $X_r$ is smooth, each $X_i$ is a toric variety, and the map $X_{i + 1} \rightarrow X_i$ is the blowup along a smooth $T$-invariant subvariety of $X_i$.  We proceed by induction on $r$ and the dimension of $X$.

Let $X' = X(\Delta') = X_1$.  By induction on $r$, we may assume that $\bigoplus_{\sigma' \in \Delta'} \iota_{\sigma'}^*$ maps $A^*_T(X')$ isomorphically onto $\PP^*(\Delta')$.  Since $\pi$ maps $O_{\sigma'}$ isomorphically onto $O(\sigma)$ if $\sigma'$ is a maximal cone in the subdivision of $\sigma$ induced by $\Delta'$, it follows that $\bigoplus_{\sigma \in \Delta} \iota_\sigma^*$ maps $A^*_T(X)$ injectively into $\PP^*(\Delta)$.  It remains to show that every integral piecewise polynomial function on $\Delta$ is in the image of $A^*_T(X)$.

Say $\pi$ is the blowup along $V(\tau) \subset X$, and $V(\rho) = \pi^{-1} (V(\tau))$.  Let $\Star \tau$ be the set of cones in $\Delta$ containing $\tau$, and let $\Delta_\tau$ be the fan whose cones are the projections of cones in $\Star \tau$ to $(N / N_\tau)_\R$, where $N_\tau$ is the sublattice generated by $\tau \cap N$.  Then $V(\tau)$ is the toric variety associated to $\Delta_\tau$ \cite[Section 3.1]{Fulton}.  By induction on dimension, we may assume $A^*_{T_\tau}(V(\tau)) \cong \PP^*(\Delta_\tau)$, where $T_\tau$ is the dense torus in 
$V(\tau)$.  Choosing a splitting $T \cong T_\tau \oplus T'$, we have
\[
A^*_T(V(\tau)) \ \cong \ A^*_{T_\tau}(V(\tau)) \otimes A^*_{T'}(\mathrm{pt})  \ \cong \ \PP^*(\Star \tau).
\]
Similarly, we have $A^*_T(V(\rho)) \cong \PP^*(\Star \rho)$.

Note that $\Star \rho$ is a subdivision of $\Star \tau$, and $\Delta$ and $\Delta'$ coincide away from $\Star \tau$ and $\Star \rho$.  By Lemma \ref{equivariant Kimura}, a class in $A^*_T(X')$ is in the image of $A^*_T(X)$ if and only if  its restriction to $V(\rho)$ is in the image of $A^*_T(V(\tau))$.  Therefore, a piecewise polynomial function on $\Delta'$ is in the image of $A^*_T(X)$ if and only if its restriction to $\Star \rho$ is the pullback of a piecewise polynomial function on $\Star \tau$.  In particular, the pullback of any piecewise polynomial function on $\Delta$ is in the image of $A^*_T(X)$, as required.  \hfill $\Box$

\vspace{10 pt}

\noindent \emph{Proof of Theorem \ref{GKM}:} Since $X$ is complete, $\Delta$ is a complete fan and all of its maximal cones are $n$-dimensional.  Since a piecewise polynomial function is uniquely determined by its restriction to the maximal cones, and since the orbits $O_\sigma$ for $n$-dimensional cones $\sigma$ are precisely the $T$-fixed points of $X$, the injectivity in Theorem \ref{main} gives injectivity in Theorem \ref{GKM}.  Furthermore, a collection of polynomials $f_1, \ldots, f_r$ on the maximal cones $\sigma_1, \ldots, \sigma_r$ of $\Delta$ give a piecewise polynomial function if and only if, for any codimension 1 cone $\tau$ that is the intersection of the maximal cones $\sigma_i$ and $\sigma_j$, $f_i|_\tau = f_j|_\tau$.  Since the orbits $O_\tau$ for codimension 1 cones $\tau$ are precisely the 1-dimensional orbits in $X$, and $O_{\sigma_i}$ and $O_{\sigma_j}$ are the fixed points in the closure of $O_\tau$, the fact that $A^*_T(X)$ surjects onto $\PP^*(\Delta)$ gives the second part of Theorem \ref{GKM}. \hfill $\Box$

\begin{remark}\emph{ Given that equivariant Chow cohomology classes on toric varieties are naturally identified with piecewise polynomial functions, it seems potentially interesting to ask what Chow cohomology class is represented whenever a piecewise polynomial function occurs naturally.  For instance, asymptotic cohomological functions of toric divisors give rational piecewise polynomial functions on the divisor class group of a toric variety, including the volume function, which is piecewise polynomial with respect to the Gelfand-Kapranov-Zelevinsky decomposition of the effective cone \cite{HKP}.  We do not know what cohomology classes these functions represent. }
\end{remark}

\section{Equivariant Chern classes}

This project was originally motivated by a desire to have a simple, combinatorial formula for the equivariant Chern classes of equivariant vector bundles on toric varieties. We refer the reader to \cite{Kaneyama}, \cite{Klyachko}, \cite{Knutson-Sharpe}, and \cite{Perling} for the basic facts about equivariant vector bundles on toric varieties, and simply recall that an equivariant vector bundle
\[
p: E \rightarrow X(\Delta)
\]
 is a vector bundle together with an action of the dense torus $T \subset X$ on $E$ such that $p$ is $T$-equivariant.  An equivariant vector bundle on an affine toric variety $U_\sigma$ is equivariantly trivializable, and is uniquely determined up to isomorphism by the action of the stablizer of a point $x_\sigma$ in the minimal orbit $O_\sigma$ on the fiber $E_{x_\sigma}$.  Hence a vector bundle on $U_\sigma$ is given by a multiset ${\bf u}_\sigma \subset M_\sigma$, where the multiplicity of $u$ in ${\bf u}_\sigma$ is the dimension of the $\chi^u$-isotypical component of $E_ {x_\sigma}$.
Let $E \rightarrow X(\Delta)$ be an equivariant vector bundle.  For $\sigma \in \Delta$, let ${\bf u}_\sigma \subset M_\sigma$ be the multiset determined by $E|_{U_\sigma}$. 

\begin{theorem}\label{chern classes}
The equivariant Chern class $c_i^T(E)$ in $A^i_T(X(\Delta))$ is naturally identified with the piecewise polynomial function whose restriction to $\sigma$ is the $i$-th elementary symmetric function $e_i({\bf u}_\sigma)$.
\end{theorem}

\noindent \emph{Proof:}  By Theorem \ref{main}, $c_i^T(E)$ is determined by $c_i^T(E|_{O_\sigma})$ for $\sigma \in \Delta$, and $E|_{O_\sigma}$ is equivariantly isomorphic to a sum of line bundles whose equivariant first Chern classes are given by the multiset ${\bf u}_\sigma$.

\vspace{ 10 pt }

\noindent Klyachko gives another formula for the Chern class of an equivariant vector bundle in the ordinary Chow ring $A^*(X)$, when $X$ is smooth and complete, using a resolution by sums of line bundles \cite[Section~3.2]{Klyachko}.  This formula may be interpreted as a formula in the Stanley-Reisner ring, by equating the class of a prime $T$-invariant divisor with the corresponding monomial generator, and can be deduced from Theorem \ref{chern classes} using the isomorphism $\PP^*(\Delta) \cong SR(\Delta)$ and the projection $A^*_T(X) \rightarrow A^*(X)$.

\section{Cohomology of toric prevarieties}

The arguments of this paper, suitably interpreted, work without change to compute equivariant Chow cohomology  and equivariant Chern classes on toric prevarieties in terms of piecewise polynomials.   Toric prevarieties have appeared as universal embedding spaces in W\l odarczyk's work \cite{Wlodarczyk}, and were given a combinatorial treatment in \cite{AH}.  Recall that a toric prevariety is a normal, but not necessarily separated, prevariety $X$ with a dense torus $T \subset X$ and an action of $T$ on $X$ extending the action of $T$ on itself.

Let $X$ be a toric prevariety, and $\Delta$ the set of $T$-invariant affine open subsets of $X$, partially ordered by inclusion.  For $U \in \Delta$, the set of one-parameter subgroups $k^* \rightarrow T$ that extend to regular morphisms $\A^1 \rightarrow U$ are exactly the lattice points in some convex rational polyhedral cone $\sigma \subset N_\R$.  So we have a natural map of posets
\[
\phi: \Delta \rightarrow \mathrm{Cone\, } N_\R
\]
to the poset of convex rational polyhedral cones in $N_\R$, and $\phi$ maps the set of $T$-invariant affine subsets of $U$ bijectively onto the faces of $\phi(U)$.  This exactly means that $(\Delta, \phi)$ is a (unweighted) multifan, in the sense of \cite{HM}.  When no confusion seems possible, we omit $\phi$, and write simply $\Delta$ to denote the multifan $(\Delta, \phi)$.  Note that, given a multifan $\Delta$, we can construct a toric prevariety $X(\Delta)$ by gluing together the affine toric varieties $U_{\phi(\sigma)}$, for $\sigma \in \Delta$, according to the poset data.  We say that the multifan $\Delta$ is simplicial if $\phi(\sigma)$ is simplicial for every $\sigma \in \Delta$, so $X(\Delta)$ has at worst quotient singularities if and only if $\Delta$ is simplicial.

Although the elements of $\Delta$ are, by definition, affine varieties, we use notation suggestive of cones in fans, denoting elements of $\Delta$ by $\sigma$ and $\tau$.   We define the ring of integral piecewise polynomials on the multifan $\Delta$ by
\[
  \PP^*(\Delta) = \{ f_\sigma \in \bigoplus_{\sigma \in \Delta} \Sym M_{\phi(\sigma)} : f_\sigma|_{\phi(\tau)} = f_\tau \mbox{ for } \tau \prec \sigma \}.
  \]
It is perhaps helpful to think of defining a rational polyhedral cone complex with integral structure, in the sense of \cite[pp.\ 69--70]{KKMS}, by gluing the cones $\phi(\sigma)$, for $\sigma \in \Delta$, according to the poset data.  Elements of $\PP^*(\Delta)$ can then be viewed as honest, real-valued functions on the underlying topological space of this cone complex that are given by an integral polynomial on each cone.

We write $O_\sigma$ for the minimal orbit in the affine open subvariety corresponding to $\sigma$, and $\iota_\sigma$ for the inclusion of $O_\sigma$ in $X(\Delta)$.

\begin{theorem}\label{prevarieties}
Let $X = X(\Delta)$ be the toric prevariety corresponding to a multifan~$\Delta$.  Then $\bigoplus_{\sigma \in \Delta} \iota_\sigma^*$ maps $A^*_T(X)$ isomorphically onto $\PP^*(\Delta)$.
\end{theorem}

\noindent The proof of Theorem \ref{prevarieties} is identical to the proof of Theorem \ref{main}, given that capping with the equivariant fundamental class gives an isomorphism $A^*_T(X) \cong A_*^T(X)$ for smooth toric prevarieties $X$, as we now prove.  This isomorphism does not follow immediately from \cite[Corollary 17.4]{IT}, because $X$ may not be separated.  It is known that $A^*(Y)_\Q$ is isomorphic to $A_*(Y)_\Q$ for any prevariety $Y$ with at worst quotient singularities.  However, this is a deep result and it is apparently not known whether it holds with integer coefficients when $Y$ is smooth; see \cite[Section 6]{Vistoli}, and also \cite{Kimura thesis}.  The required isomorphism is easy to construct in our case.

\begin{lemma}\label{cap isomorphism}
Let $X = X(\Delta)$ be a smooth toric prevariety.  Then capping with the fundamental class gives an isomorphism $A_T^*(X) \cong A_*^T(X)$.
\end{lemma}

\noindent \emph{Proof:} As in the Cox construction of separated toric varieties \cite{Cox}, we can construct $X$ as a geometric quotient
\[
 X \cong (\A^n \smallsetminus Z) / T',
\]
where $T^n$ is the $n$-dimensional torus acting coordinatewise on $\A^n$, and $T' \subset T^n$ acts freely on the complement of the coordinate subspace arrangement $Z$.  Hence, $A^*_T(X)$ is canonically identified with $A^*_{T^n}(\A^n \smallsetminus Z)$.  Since $\A^n \smallsetminus Z$ is smooth and separated, capping with the equivariant fundamental class gives $A^*_{T^n}(\A^n \smallsetminus Z) \cong A_*^{T^n}(\A^n \smallsetminus Z)$.  Since $A_*^{T^n}(\A^n \smallsetminus Z)$ is canonically identified with $A_*^T(X)$, the lemma follows. \hfill $\Box$

\vspace{10 pt}

\noindent The equivariant Chern classes of equivariant vector bundles on toric prevarieties are also given as in Theorem~\ref{chern classes}.

\section{Cohomology of hypertoric varieties}

We now apply Theorem \ref{prevarieties} to compute the equivariant Chow cohomology of general complex hypertoric varieties.  We briefly recall the construction of general hypertoric varieties (by ``general" we mean that we consider only regular values of the moment map), and refer the reader to \cite{BD}, \cite{HS}, and \cite{PW} for details and further references, as well as for the construction of hypertoric varieties associated to nonregular values of the moment map.  Although most of the following construction works over an arbitrary field, we will eventually use the cycle class map in an essential way in the proof of Theorem \ref{hypertoric}.  For simplicity, we restrict to working over the complex numbers for the remainder of this paper.  Let $T^n = (\C^*)^n$ be the $n$-dimensional complex torus, acting coordinatewise on $\C^n$.  
                                                                                                                                                                                                                                                                                                                                                                                                                                                                                                                                                                                                                                                                                                                                                                                                                                                                                                                                                                                                                                                                                                                                                                                                                                                                                                                                                           
Let $T' \subset T^n$ be a subtorus, and let $T = T^n / T'$ be the quotient torus.  Let 
\[
  0 \rightarrow N' \rightarrow \Z^n \rightarrow N \rightarrow 0
\]
be the induced short exact sequence of lattices of one parameter subgroups, with
\[
0 \leftarrow M' \leftarrow \Z^n \leftarrow M \leftarrow 0
\]
the dual short exact sequence of character lattices.
For $i = 1, \ldots, n$, let $u_i$ be the image of $e_i$ in $M'$, let $v_i$ be the image of $e_i$ in $N$, and let $\cA$ be the cooriented hyperplane arrangement in $M_\R$ dual to $\{v_1, \ldots, v_n \}$.

The big torus $T^n$ acts naturally on the cotangent bundle $T^* \C^n \cong \C^n \times \C^n$ by $t(z,w) = (tz, t^{-1} w)$.  The restriction of this action to $T'$ induces a moment map
\[
\mu: \C^n \times \C^n \rightarrow M'_\C
\]
given by 
\[
(z,w) \mapsto \sum_{i = 1}^n (z_i w_i) \cdot u_i.
\]
The regular values for $\mu$ are exactly those $\lambda$ in $M'_\C$ that are not contained in any of the hyperplanes spanned by subsets of $\{ u_1, \ldots, u_n \}$.  For such $\lambda$, $T'$ acts locally freely on $\mu^{-1}(\lambda)$, hence we may define the hypertoric variety $\M_\lambda(\cA)$ to be the geometric quotient
\[
  \M_\lambda(\cA) = \mu^{-1}(\lambda) /\, T'.
\]
Since $\mu$ commutes with the action of $T^n$, $\M_\lambda(\cA)$ inherits a natural action of the quotient torus $T$.

Proudfoot and Webster observed that the hypertoric variety $\M_\lambda(\cA)$ is naturally a $T$-equivariant affine bundle over a toric prevariety $X(\Delta_\cA)$ \cite[Remark~3.6]{PW}, as we now describe.  Let $\Delta_\cA$ be the set of linearly independent subsets of $\{v_1, \ldots, v_n\}$, partially ordered by inclusion.  We have a natural map of posets
\[
\Delta_\cA \rightarrow \mathrm{Cone\, }N_\R,
\]
taking a collection of linearly independent vectors to the cone that they span in $N_\R$, making $\Delta_\cA$ a multifan.  Now $X(\Delta_\cA)$ can be constructed as a geometric quotient, as in the Cox construction of simplicial toric varieties \cite{Cox},
\[
  X(\Delta_\cA) = (\C^n \smallsetminus Z) /\, T',
\]  
where $Z$ is the subspace arrangement cut out by the monomials $z_{i_1}\! \cdots z_{i_r}$ such that $\{v_{i_1}, \ldots, v_{i_r}\}$ does not span a cone in $\Delta_\cA$. (As noted by Proudfoot, the prevariety $X(\Delta_\cA)$ may also be constructed as the union of all of the GIT quotients $\C^n /\!/_\alpha\, T'$ for different generic linearizations $\alpha$ of $T'$ acting on the trivial line bundle over $\C^n$, glued together along the open subsets where they agree.)  It is now straightforward to check, using the fact that  $\{v_{i_1}, \ldots, v_{i_r}\}$ spans a cone in $\Delta_\cA$ if and only if it is linearly independent, that projection to the first factor $p_1 : \C^n \times \C^n \rightarrow \C^n$ maps $\mu^{-1}(\lambda)$ equivariantly onto $\C^n \smallsetminus Z$, inducing an equivariant affine bundle $M_\lambda(\cA) \rightarrow X(\Delta_\cA)$ that trivializes equivariantly over each $T$-invariant affine open $U_\sigma$ in $X(\Delta_\cA)$.  

\begin{remark} \emph{ The above construction shows that any simplicial toric prevariety is the base of an affine bundle with separated total space: suppose $\Delta$ is a simplicial multifan and $\{ v_1, \ldots, v_r \}$ is the multiset of primitive generators of rays of $\Delta$, with $\cA$ the dual cooriented hyperplane arrangement.  Then $\Delta$ is a submultifan of $\Delta_\cA$, inducing an inclusion $X(\Delta) \subset X(\Delta_\cA)$.  If $p: \M_\lambda(\cA) \rightarrow X(\Delta_\cA)$ is the affine bundle associated to a regular value $\lambda$, then $p^{-1} \left(X(\Delta)\right)$ is a separated affine bundle over $X(\Delta)$. }
\end{remark}

It is apparently not known in general whether affine bundles induce isomorphisms in Chow cohomology.    However, the methods of this paper yield this result for the affine bundles $\M_\lambda(\cA) \rightarrow X(\Delta_\cA)$.

\begin{theorem}\label{hypertoric}
Let $p: \M \rightarrow X$ be an equivariant affine bundle with separated total space over a complex toric prevariety.  Then $p^*: A^*_T(X) \rightarrow A^*_T(\M)$ is an isomorphism.
\end{theorem}

\noindent \textbf{Corollary}
\emph{For any regular value $\lambda$, the equivariant Chow cohomology ring of the complex hypertoric variety $\M_\lambda(\cA)$ is naturally isomorphic to $\PP^*(\Delta_\cA)$.}

\vspace{10 pt}

In preparation for the use of the cycle class map in the proof of Theorem \ref{hypertoric}, we first prove a lemma on the equivariant singular cohomology of complex toric prevarieties and affine bundles.

\begin{lemma}\label{singular lemma}
Let $X = X(\Delta)$ be a smooth complex toric prevariety, and let $p : \M \rightarrow X$ be an equivariant affine bundle.  Then $\bigoplus_{\sigma \in \Delta} \iota_\sigma^*$ maps $H^*_T(X)$ isomorphically onto $\PP^*(\Delta)$, and $p^*$ maps $H^*_T(X)$ isomorphically onto $H^*_T(\M)$.  

If $\Delta$ is simplicial, but $X$ is not necessarily smooth, the analogous results hold with $\Q$ coefficients.
\end{lemma}

\noindent \emph{Proof:}  The isomorphism $H^*_T(X) \cong \PP^*(\Delta)$ follows from the degeneration of the Mayer-Vietoris spectral sequence for the cover of $X \times_T ET$ by $\{U_\sigma \times_T ET\}_{\sigma \in \Delta}$, just as for separated toric varieties.  The fact that $p^*$ is an isomorphism follows from the Leray-Serre spectral sequence applied to the affine fiber bundle
\[
\M \times_T ET \rightarrow X \times_T ET.  
\]
The proof for $\Delta$ simplicial, with rational coefficients, is identical.
\hfill $\Box$

\vspace{10 pt}

For a cooriented hyperplane arrangement $\cA$ as above, the standard argument relating piecewise polynomials to Stanley-Reisner rings shows that $\PP^*(\Delta_\cA)_\Q$ is naturally isomorphic to $SR(\Delta_\cA)_\Q$.  So Lemma~\ref{singular lemma}, applied to the affine bundle $\M_\lambda(\cA) \rightarrow X(\Delta_\cA)$, gives another proof of the fact that the equivariant singular cohomology ring $H^*_T(\M_\lambda(\cA))_\Q$ is naturally isomorphic to $SR(\Delta_\cA)_\Q$ (and with integer coefficients if $\M_\lambda(\cA)$ is smooth).  This statement appears in \cite[Theorem 6.1]{PW}, and is due to Hausel and Sturmfels for the general case (it follows from arguments similar to the ordinary cohomology computations in \cite{HS}), and to Konno for the smooth case with integer coefficients \cite{Konno}.

\vspace{10 pt}

\noindent \emph{Proof of Theorem \ref{hypertoric}:}  We first consider the case when $X = X(\Delta)$ is smooth.  Since $\M$ is separated and smooth, capping with $[\M]_T$ gives an isomorphism $A^*_T(\M) \cong A_*^T(\M)$.  By Lemma \ref{cap isomorphism}, capping with $[X]_T$ gives $A^*_T(X) \cong A_*^T(X)$, as well.  The flat pullback map $A_*^T(X) \rightarrow A_*^T(\M)$ is surjective (apply \cite[Proposition 1.9]{IT} to the affine bundles $\M_T \rightarrow X_T$), so the composition
\[
A^*_T(X) \cong A_*^T(X) \rightarrow A_*^T(\M) \cong A^*_T(\M)
\]
is surjective.  We must show that $p^*$ is injective.  Consider the cycle class map
\[
\mathrm{cl} : A^*_T(\M) \rightarrow H^*_T(\M) \cong \PP^*(\Delta),
\]
where the latter isomorphism is given by Lemma \ref{singular lemma}.  The composition $\mathrm{cl} \circ p^*$ is the canonical isomorphism $A^*_T(X) \cong \PP^*(\Delta)$, hence $p^*$ is also injective.

If $X$ is not smooth, choose a toric resolution of singularities
\[
X' = X(\Delta') \stackrel{\pi}{\longrightarrow} X,
\]
and let $\M' \rightarrow \M$ be the induced resolution of $\M$, where $\M' = \M \times_X X'$.  Since $A^*_T(X) \rightarrow A^*_T(\M') \cong \PP^*(\Delta')$ is injective, $p^*: A^*_T(X) \rightarrow A^*_T(\M)$ is injective as well.  To see that  $p^*$ is surjective, note that $\M' \rightarrow \M$, being the base change of an envelope, is an envelope, and hence $A^*_T(\M) \rightarrow A^*_T(\M') \cong \PP^*(\Delta')$ is injective.  Restricting cohomology classes to $p^{-1} O_\sigma \subset \M$ for $\sigma \in \Delta$ shows that the image of $A^*_T(\M)$ is contained in $\PP^*(\Delta)$.  Since $A^*_T(X)$ surjects onto $\PP^*(\Delta)$, it follows that $p^*$ is surjective, as required.  \hfill $\Box$

\vspace{10 pt}

\noindent \textbf{Acknowledgments}\, I thank Nick Proudfoot for stimulating my interest in toric prevarieties and hypertoric varieties, and for many helpful discussions.  I am grateful to the organizers of GAEL XIII, where much of this research was done, for their hospitality.  This work was supported by a Graduate Research Fellowship from the NSF.

\vspace{10 pt}

\noindent University of Michigan \\
\noindent Department of Mathematics \\
\noindent 2074 East Hall, 530 Church St. \\
\noindent Ann Arbor, MI 48109

\vspace{5 pt}

\noindent \emph{sdpayne@umich.edu}

\end{document}